\documentclass{elsart3-1}

\usepackage{amssymb,amsmath}

\usepackage[english,francais]{babel}

\newtheorem{theorem}{Theorem}[section]
\newtheorem{lemma}[theorem]{Lemma}
\newtheorem{e-proposition}[theorem]{Proposition}

\newtheorem{e-definition}[theorem]{Definition\rm}

\renewcommand{\theequation}{\arabic{equation}}
\renewcommand{\epsilon}{\varepsilon}
\def\og{\leavevmode\raise.3ex\hbox{$\scriptscriptstyle\langle\!\langle$~}}
\def\fg{\leavevmode\raise.3ex\hbox{~$\!\scriptscriptstyle\,\rangle\!\rangle$}}

\begin{document}

\begin{frontmatter}


\selectlanguage{english}
\title{\vspace*{-30mm}Stable anti-Yetter-Drinfeld modules}




\selectlanguage{english}
\author[authorlabel1,authorlabel2]{Piotr M.~Hajac}
\ead{http://www.fuw.edu.pl/$\!\widetilde{\phantom{m}}\!$pmh}
\author[authorlabel3]{Masoud Khalkhali}
\ead{masoud@uwo.ca}
\author[authorlabel3]{Bahram Rangipour}
\ead{brangipo@uwo.ca}
\author[authorlabel4]{Yorck Sommerh\"auser}
\ead{sommerh@mathematik.uni-muenchen.de}

\address[authorlabel1]{Instytut Matematyczny, Polska Akademia Nauk,
ul.\ \'Sniadeckich 8, Warszawa, 00-956 Poland}
\address[authorlabel2]{Katedra Metod Matematycznych Fizyki, Uniwersytet Warszawski,
ul.\ Ho\.za 74, Warszawa, 00-682 Poland}
\address[authorlabel3]{Department of Mathematics, University of
Western Ontario, London ON, Canada}
\address[authorlabel4]{Mathematisches Institut, Universit\"at M\"unchen,
Theresienstr.\ 39, 80333 M\"unchen,  Germany}

\begin{abstract}
We define and study a class of entwined modules (stable anti-Yetter-Drinfeld modules)
that serve as coefficients for
the Hopf-cyclic homology and cohomology. 
In particular, we explain their relationship with Yetter-Drinfeld modules and Drinfeld doubles.
Among sources of examples of stable anti-Yetter-Drinfeld modules, we find Hopf-Galois
extensions with a flipped version of the Miyashita-Ulbrich action. 

{\it To cite this article: P.~M.~Hajac et al., C. R.
Acad. Sci. Paris, Ser. I 336 (2003).}

\vskip 0.5\baselineskip

\selectlanguage{francais}
\noindent{\bf R\'esum\'e}
\vskip 0.5\baselineskip
\noindent
{\bf  Modules anti-Yetter-Drinfeld stables}
 Nous d\'efinissons et \'etudions une classe de modules enlac\'es
(modules anti-Yetter-Drinfeld stables) qui servent de coefficients
pour l'homologie et la cohomologie Hopf-cyclique.
En
particulier, nous expliquons leurs liens avec les modules de
Yetter-Drinfeld et les doublets de Drinfeld. Parmi les sources d'exemples
de modules anti-Yetter-Drinfeld stables, nous trouvons des
extensions de Hopf-Galois munies d'une version transpos\'ee de l'action de
Miyashita-Ulbrich.

 {\it Pour citer cet article~: P.~M.~Hajac et al., C. R.
Acad. Sci. Paris, Ser. I 336 (2003).}

\end{abstract}
\end{frontmatter}

\selectlanguage{english}
\section{Introduction}
 The aim of this paper is to define and provide sources of examples
of stable anti-Yetter-Drinfeld modules. They play the role of coefficients for Hopf-cyclic 
theory~\cite{hkrsb}.
In particular, we claim that modular pairs in involution of Connes and Moscovici are precisely 1-dimensional
stable anti-Yetter-Drinfeld modules.

Throughout the paper we assume that $H$
is a Hopf algebra with a bijective antipode.
 On the one hand, the bijectivity of the antipode is implied
by the existence of a modular pair in involution, so that then it need not be
assumed. On the other hand, some parts of arguments might work even if the
antipode is not bijective. We avoid such discussions.
The  coproduct, counit and antipode of $H$ are denoted by
$\Delta$,  $\epsilon$ and $S$,  respectively.
For the coproduct we use the  notation $\Delta(h)=h^{(1)}\otimes h^{(2)}$,
for a left coaction on $M$ we write $_M\Delta(m)=m^{(-1)}\otimes m^{(0)}$, and for
a right coaction  $\Delta_M(m)=m^{(0)}\otimes m^{(1)}$.  The summation
 symbol is suppressed everywhere. We assume all algebras to be associative, unital
and over the same ground field $k$. The symbol ${\mathcal O}(X)$ stands for the
algebra of polynomial functions  on $X$.

P.M.H.\ acknowledges  the Marie Curie Fellowship
 HPMF-CT-2000-00523 and KBN grant 2 P03A 013 24.
All  authors are grateful to T.~Brzezi\'nski and M.~Furuta
for their help and comments, and to D.~Perrot for the French translation.


\section{The transformation of Yetter-Drinfeld modules}

It turns out that, in order to incorporate coefficients into  cyclic theory,
we  need to alter the concept of a Yetter-Drinfeld module by replacing
  the antipode by its inverse in the Yetter-Drinfeld compatibility condition
 between actions and coactions. We call the modules-comodules satisfying
 the thus modified  Yetter-Drinfeld compatibility condition {\em anti-Yetter-Drinfeld
 modules}\footnote{
 This concept was  devised independently by Ch.~Voigt
and, also independently, by P.~Jara and D.~\c{S}tefan.}.
   Just as Yetter-Drinfeld modules
 come in 4 different versions depending on the side of actions and coactions
 (see \cite[p.181]{cmz02} for a general formulation),
 so do the anti-Yetter-Drinfeld modules.
All versions are completely equivalent and can be derived from one another
 by replacing a Hopf algebra
 $H$ by $H^{cop}$, $H^{op}$, or $H^{op,cop}$, respectively.\\
\vspace*{-2mm}

\begin{e-definition}
\label{ayd}
Let $H$ be a Hopf algebra with a bijective antipode $S$,
 and $M$ a module and comodule
 over $H$. We call $M$ an anti-Yetter-Drinfeld module iff
 the action and coaction are compatible in the following sense:
\begin{align}
&\label{ayd1}
_M\Delta(hm)=h^{(1)}m^{(-1)}S^{-1}(h^{(3)})\otimes h^{(2)}m^{(0)}&
&\mbox{\em if $M$ is a left module and a left comodule},&
\\ &
\Delta_M(hm)=h^{(2)}m^{(0)}\otimes h^{(3)}m^{(1)}S(h^{(1)})&
&\mbox{\em if $M$ is a left module and a right comodule},&
\\ &
_M\Delta(mh)=S(h^{(3)})m^{(-1)}h^{(1)}\otimes m^{(0)}h^{(2)}&
&\mbox{\em if $M$ is a right module and a left comodule},&
\\ &\label{ayd4}
\Delta_M(mh)= m^{(0)}h^{(2)}\otimes S^{-1}(h^{(1)})m^{(1)}h^{(3)}&
&\mbox{\em if $M$ is a right module and a right comodule}.&
\end{align}
\end{e-definition}~\vspace*{-2mm}

To make  cyclic theory work, we also need to assume that the action splits
coaction, i.e., for all $m\in M$,
$m^{(-1)}m^{(0)}=m$, $m^{(1)}m^{(0)}=m$, $m^{(0)}m^{(-1)}=m$,
$m^{(0)}m^{(1)}=m$, for the left-left, left-right, right-left, and right-right
versions, respectively. We call modules satisfying this condition {\em stable}.
Let us emphasize that it is the anti-Yetter-Drinfeld condition rather than
the Yetter-Drinfeld condition that makes the homomorphism
$action\circ coaction$ $H$-linear and $H$-colinear. Therefore the stability condition
$action\circ coaction= {\rm id}$ suits the former and not the latter.
The first class of examples of stable anti-Yetter-Drinfeld\ modules is provided by modular
pairs in involution \cite[p.8]{cm00}. Since such pairs occur naturally in different contexts,
Lemma~\ref{1dim} and Lemma~\ref{affine}
 guarantee ample amount of examples of
anti-Yetter-Drinfeld\ modules.\\
\vspace*{-2mm}

\begin{lemma}\label{1dim}
Let the ground field
$k$ be a right module over $H$ via a character $\delta$ and a left comodule
over $H$ via a group-like $\sigma$. Then $k=^\sigma\!\!\!k_\delta$ is a stable anti-Yetter-Drinfeld\ module
{\em if and only if} $(\delta,\sigma)$ is a modular pair in involution.
\end{lemma}~\vspace*{-2mm}

The anti-Yetter-Drinfeld modules do not form a monoidal category themselves, but rather a so-called $\mathcal C$-category over the category of Yetter-Drinfeld modules (see \cite[p.351]{p-b77} for details). More precisely:\\
\vspace*{-2mm}

\begin{lemma}\label{affine}
Let $N$ be a Yetter-Drinfeld module and $M$ an anti-Yetter-Drinfeld~ module.
Then $N\otimes M$  is an anti-Yetter-Drinfeld~ module via
$h(n\otimes m)=h^{(1)}n\otimes h^{(2)}m$,
$_{N\otimes M}\Delta(n\otimes m)=n^{(-1)}m^{(-1)}\otimes n^{(0)}\otimes m^{(0)}$,
for the left-left case, and via
$h(n\otimes m)=h^{(2)}n\otimes h^{(1)}m$,
$\Delta_{N\otimes M}(n\otimes m)=n^{(0)}\otimes m^{(0)}\otimes n^{(1)}m^{(1)}$,
for the left-right case. Similarly, $M\otimes N$ is an anti-Yetter-Drinfeld\ module via
$(m\otimes n)h=mh^{(2)}\otimes nh^{(1)}$,
$_{M\otimes N}\Delta(n\otimes m)=m^{(-1)}n^{(-1)}\otimes m^{(0)}\otimes n^{(0)}$,
for the right-left case, and via
$(m\otimes n)h=mh^{(1)}\otimes nh^{(2)}$,
$\Delta_{M\otimes N}(m\otimes n)=m^{(0)}\otimes n^{(0)}\otimes m^{(1)}n^{(1)}$,
for the left-right case.
\end{lemma}~\vspace*{-2mm}

Note that, just as the right-right Yetter-Drinfeld modules are entwined modules \cite{b-t99}
for the entwining $\psi(h'\otimes h)=h^{(2)}\otimes S(h^{(1)})h'h^{(3)}$, the right-right
anti-Yetter-Drinfeld modules are entwined
with respect to $\psi(h'\otimes h)=h^{(2)}\otimes S^{-1}(h^{(1)})h'h^{(3)}$. (Other cases are completely analogous.)

An intermediate step between modular pairs in involution and stable anti-Yetter-Drinfeld modules
is given by matched and comatched pairs of \cite{kr03}.
Whenever the antipode is equal to its inverse, the difference between the
Yetter-Drinfeld and anti-Yetter-Drinfeld\ conditions disappears.
For a group ring Hopf algebra $kG$, a left
$H$-comodule
 is simply a $G$-graded vector space $M={\bigoplus}_{g\in G}M_g$,
  where the coaction is defined by $M_g\ni m\mapsto g\otimes m$.
 An action of $G$ on $M$ defines an (anti-)Yetter-Drinfeld
 module if and only if
for all $g,h\in G$ and $m\in M_g$ we have
 $hm\in M_{hgh^{-1}}$. The stability condition means simply that $gm=m$ for all
$g\in G$, $m\in M_g$.
A very concrete classical example of a stable (anti-)Yetter-Drinfeld module is
provided by the Hopf fibration. Then $H={\mathcal O}(SU(2))$ and $M={\mathcal O}(S^2)$. Since
$S^2\cong SU(2)/U(1)$, we have a natural left action of $SU(2)$ on $S^2$.
Its pull-back
makes $M$ a left $H$-comodule. On the other hand, one can view $S^2$ as the set
of all traceless matrices of $SU(2)$. The pull-back of this embedding
$j:S^2{\hookrightarrow}SU(2)$ together with the multiplication in ${\mathcal O}(S^2)$
defines a left action of $H$ on $M$. It turns out that the equivariance property
$j(gx)=gj(x)g^{-1}$ guarantees the anti-Yetter-Drinfeld\ condition, and this combined with
the injectivity of $j$ ensures the stability of $M$. This stability mechanism
can be generalized in the following way.\\
\vspace*{-2mm}

\begin{lemma}
 Let $M$ be an algebra and a left $H$-comodule. Assume that
   $\pi: H\rightarrow M$ is an epimorphism of algebras and
   the action $hm=\pi(h)m$ makes $M$ an anti-Yetter-Drinfeld\ module. Assume
also that $\pi(1^{(-1)})1^{(0)}=1$. Then
   $M$ is a stable module.
\end{lemma}

\section{Hopf-Galois extensions and the opposite Miyashita-Ulbrich action}

Another source of examples is provided by Hopf-Galois theory. These examples are
purely quantum in the sense that the employed actions are automatically trivial
for commutative algebras. To fix the notation and terminology, recall that an algebra
and an $H$-comodule is called a comodule algebra if the coaction is an algebra
homomorphism. An $H$-extension $B:=\{p\in P\;|\;\Delta_P(p)=p\otimes 1\}\subseteq P$ is called
Hopf-Galois iff the canonical map $can: P\otimes_BP\rightarrow P\otimes H$, $can(p\otimes p')=
p\Delta(p')$, is bijective.
The bijectivity assumption
allows us to define the translation map $T:H\rightarrow P\otimes_BP$, $T(h):=can^{-1}(1\otimes h)
=:h^{[1]}\otimes_Bh^{[2]}$ (summation suppressed). It can be shown that when
everything is over a field (our standing assumption), the centralizer
$Z_B(P):=\{p\in P\;|\;bp=pb, \forall\, b\in B\}$ of $B$ in $P$ is a subcomodule of $P$.
On the other hand,
 the formula $ph=h^{[1]}ph^{[2]}$ defines a right action on $Z_B(P)$ called
the Miyashita-Ulbrich action. This action and coaction satisfy the
Yetter-Drinfeld compatibility condition \cite[(3.11)]{dt89}. The following proposition modifies
the Miyashita-Ulbrich action  so as to obtain stable anti-Yetter-Drinfeld modules.\\
\vspace*{-2mm}

\begin{e-proposition}
\label{mu}
Let  $B\subseteq P$ be a Hopf-Galois
$H$-extension such that $B$ is central in $P$. Then $P$ is a right-right
stable anti-Yetter-Drinfeld\ module via the action
$
ph=(S^{-1}(h))^{[2]}p(S^{-1}(h))^{[1]}
$
and the right coaction on $P$.
\end{e-proposition}~\vspace*{-2mm}

The simplest examples are obtained for $P=H$. A broader class is given by the
so-called Galois objects \cite{c-s98}. Then quantum-group coverings at roots of unity
provide examples with central coinvariants bigger than the ground field
(see \cite{dhs99} and examples therein). Finally, one can generalize Proposition~\ref{mu}
to arbitrary Hopf-Galois extensions by replacing $P$ by $P/[B,P]$ 
\cite[Remark~4.2]{js}.

\section{The Drinfeld double comodule algebra}

For finite-dimensional Hopf algebras, the Yetter-Drinfeld modules can be understood as modules over the Drinfeld double \cite[p.220]{k-c95}. Much in the same way, 
the anti-Yetter-Drinfeld modules can also be understood as modules over a certain algebra. This makes the usual notions and operations for modules, like projectivity or induction, directly available for anti-Yetter-Drinfeld modules.
To this end,
the comodule structure of an anti-Yetter-Drinfeld module has to be converted into a module structure over the dual Hopf algebra~$H^*$, so that  from now on we assume that the Hopf algebra~$H$ is finite-dimensional. \\
\vspace*{-2mm}

\begin{e-proposition}
Let $H$ be a finite-dimensional Hopf algebra. The formula
\begin{equation}
(\varphi \otimes h)(\varphi' \otimes h') =
\varphi'{}^{(1)}(S^{-1}(h^{(3)})) \varphi'{}^{(3)}(S^2(h^{(1)})) \; 
\varphi \varphi'{}^{(2)} \otimes h^{(2)} h'
\end{equation}
turns the vector space $A(H):= H^* \otimes H$ 
into an associative algebra with the unit $\varepsilon \otimes 1$.
\end{e-proposition}~\vspace*{-2mm}

Note that the above product differs from the product in the Drinfeld double of~$H$ \cite[p.214]{k-c95} only by the additional squared antipode in the second factor. To
relate the modules over~$A(H)$ with anti-Yetter-Drinfeld modules, recall first that every right $H$-comodule~$M$ becomes a left $H^*$-module via
$\varphi m := \varphi(m^{(1)}) m^{(0)}$.
Conversely, any left $H^*$-module yields a right $H$-comodule via
$\Delta_M(m) = \sum_{i=1}^n h_i^* m \otimes h_i$. Here
$\{h_1,\ldots,h_n\}$ is a basis of~$H$ and $\{h_1^*,\ldots,h_n^*\}$ is the dual basis. 
(Of course, this comodule structure does not depend on the choice of a basis.) Using this, we get the following connection between the modules over~$A(H)$ and anti-Yetter-Drinfeld modules:\\
\vspace*{-2mm}

\begin{e-proposition}
Let $H$ be a finite-dimensional Hopf algebra.
If $M$ is a left-right anti-Yetter-Drinfeld module, it becomes a left $A(H)$-module by
$(\varphi \otimes h)m := \varphi((hm)^{(1)}) \, (hm)^{(0)}$.
Conversely, 
if $M$ is a left $A(H)$-module, it becomes a left-right anti-Yetter-Drinfeld module by
$hm := (\varepsilon \otimes h)m$,
$\Delta_M(m) := \linebreak\sum_{i=1}^n (h_i^* \otimes 1) m \otimes h_i$.
Here $\{h_1,\ldots,h_n\}$ is a basis of~$H$ and $\{h_1^*,\ldots,h_n^*\}$ its dual basis.
\end{e-proposition}~\vspace*{-2mm}

The claim of Lemma~\ref{affine} is reflected in the fact that although $A(H)$ is not a Hopf algebra
itself, it can be shown that  the formula $(\varphi \otimes h)
\mapsto (\varphi^{(2)} \otimes h^{(1)}) \otimes 
(\varphi^{(1)} \otimes h^{(2)})$ makes  $A(H)$ a right comodule algebra over the Drinfeld double~$D(H)$.

\end{document}